\documentclass[11]{amsart}
\usepackage{amssymb,amsmath,amsthm}

\usepackage{setspace}

\usepackage{lineno}

\usepackage{amssymb,amsmath}

\usepackage[english,francais]{babel}

\newtheorem{theorem}{Theorem}

\newtheorem{e-proposition}[theorem]{Proposition}

\newtheorem{e-definition}[theorem]{Definition\rm}

\providecommand{\bysame}{\leavevmode\hbox to3em{\hrulefill}\thinspace}

\def\og{\leavevmode\raise.3ex\hbox{$\scriptscriptstyle\langle\!\langle$~}}
\def\fg{\leavevmode\raise.3ex\hbox{~$\!\scriptscriptstyle\,\rangle\!\rangle$}}

\begin{document}


\selectlanguage{english}
\title{Flag structure for operators in the Cowen-Douglas class
}
\subjclass[2000]{Primary 47C15, 47B37; Secondary 47B48, 47L40}
\keywords
{The Cowen-Douglas class,  strongly irreducible, homogeneous operator, curvature, second fundamental form} 



\selectlanguage{english}
\author[Kui Ji]{Kui Ji}
\email[Kui Ji]{jikuikui@gmail.com}
\author[C. Jiang]{Chunlan Jiang}
\email[C. Jiang]{cljiang@hebtu.edu.cn}
\author[D. K. Keshari]{Dinesh Kumar Keshari}
\email[D. K. Keshari]{kesharideepak@gmail.com}
\author[G. Misra]{Gadadhar Misra}
\email[G. Misra]{gm@math.iisc.ernet.in}

\thanks{Kui Ji was supported by the Foundation for the Author of National Excellent Doctoral Dissertation of China, Grant
No. 201116. C. Jiang was supported by  National Natural Science Foundation of China, Grant No. A010602. D. K. Keshari was supported by a Research Fellowship of the Indian Institute of Science. G. Misra was supported by the J C Bose National Fellowship and the UGC, SAP -- IV}

\address[Kui Ji and Chunlan Jiang]{Department of Mathematics, Hebei Normal University, Shijiazhuang, Hebei 050016}
\address[D. K. Keshari]{Department of Mathematics\\
Texas A\& M University, College Station, TX 77843}
\address[G. Misra]{Department of Mathematics, Indian  Institute of Science,
Bangalore 560 012}




\begin{abstract}
The explicit description of homogeneous operators and localization of a Hilbert module naturally leads to the definition of a class 
of Cowen-Douglas operators possessing a flag structure. These operators are  irreducible.  We show that the flag structure is rigid in the sense that the unitary equivalence class of the operator and the flag structure determine each other.  We obtain a complete set of unitary invariants which are somewhat more tractable than those of an arbitrary operator in the Cowen-Douglas class.   


%
%

\end{abstract}


\maketitle
{\doublespacing



The Cowen-Douglas class $B_n(\Omega)$ consists of those bounded linear operators $T$ on a complex separable Hilbert space $\mathcal H$ which 
possess an open set {$\Omega\subset\mathbb C$} of eigenvalues of constant multiplicity {$n$} and  admit a holomorphic choice of eigenvectors: $s_1(w), \ldots , s_n(w)$, $w\in \Omega$, in other words, {there exists holomorphic functions $s_1, \ldots , s_n: \Omega \to \mathcal H$ which span the eigenspace of $T$ at $w\in \Omega$}.

The holomorphic choice of eigenvectors {$s_1, \ldots , s_n$} defines a holomorphic Hermitian vector bundle {$E_T$} via the map 
%
$$s: \Omega \to \mbox{\rm Gr}(n, \mathcal H),\,\, s(w) = \ker (T-w) \subseteq \mathcal H.$$
In the paper \cite{CD1},  Cowen and Douglas show that  there is a one to one correspondence between the unitary equivalence class of the operators {$T$} in $B_n(\Omega)$ and the equivelence classes of the holomorphic Hermitian vector bundles {$E_T$} determined by them. 
They also find a set of complete invariants for this equivalence consisting of the curvature {$\mathcal K$} of {$E_T$}  and a certain number of its covariant derivatives.   
Unfortunately, these invariants are not easy to compute unless {$n$ is $1$.}  
Also, it is difficult to determine, in general, when an operator {$T$} in $B_n(\Omega)$ is irreducible, again except in the case {$n=1$.} In this case, the rank of the vector bundle is {$1$} and therefore it is irreducible and so is the operator {$T$}. 

Finding similarity invariants for operators in the class $B_n(\Omega)$ has been somewhat difficult from the beginning. Counter examples to the similarity conjecture in \cite{CD1} were  given in \cite{CM1, CM2}. More recently, significant progress on the question of similarity has been made (cf. \cite{DKS, jj, jw} ). 

We isolate a subset of irreducible operators in the Cowen-Douglas class $B_n(\Omega)$ for which a complete set of tractable  unitary invariants is relatively easy to identify.  We also determine when two operators in this class are similar.  

We introduce below this smaller  class $\mathcal FB_2(\Omega)$ of operators in $B_2(\Omega)$ leaving out the more general definition of the class $\mathcal FB_n(\Omega),\, n>2,$ for now.  

\begin{e-definition}  We let $\mathcal FB_2(\Omega)$ denote the set of operators $T\in B_2(\Omega)$ which 
admit a decomposition of the form 
$$T=\begin{pmatrix}
 T_0 & S \\
 0 & T_1 \\
\end{pmatrix}$$
for some choice of operators  $T_0,T_1\in\mathcal{B}_1(\Omega)$ and  a non-zero intertwining operator $S$ between $T_0$ and $T_1$, that is, $T_0S=ST_1$.  
\end{e-definition}

An operator $T$ in  $B_2(\Omega)$ admits a decomposition of the form (cf.  \cite[Theorem 1.49, pp. 48]{jw})  $\Big (\begin{smallmatrix}
T_0 & S \\
0 & T_1 \\
\end{smallmatrix}\Big )$
for some pair of operators $T_0$ and $T_1$ in ${B}_1(\Omega)$. In defining the new class $\mathcal FB_2(\Omega)$, we are merely imposing one additional condition, namely that $T_0S=ST_1$.

We show that $T$ is in the class $\mathcal FB_2(\Omega)$ if and only if there exist a frame 
$\{\gamma_0,\gamma_1\}$ of the vector bundle $E_{T}$ such
that $\gamma_0(w)$ and $t_1(w):=\tfrac{\partial}{\partial w}\gamma_0(w)-\gamma_1(w)$ are orthogonal for all $w$ in $\Omega$.  This is also equivalent to the existence of a  frame 
$\{\gamma_0,\gamma_1\}$ of the vector bundle $E_{T}$ such
that $\tfrac{\partial}{\partial w}\|\gamma_0(w)\|^2=\langle\gamma_1(w),\gamma_0(w)\rangle,\,w\in \Omega$. Our first main theorem on unitary classification is given below, where we have set 
$\mathcal K_{T_0}(z) = -\tfrac{\partial^2}{\partial z \partial\bar{z}} \log \|\gamma_0(z)\|^2.$

\setcounter{theorem}{0}

\begin{theorem}\label{maint}
Let $T=\begin{pmatrix}T_0 & S \\
0 & T_1 \\
\end{pmatrix}$
and $\tilde{T}=\begin{pmatrix}\tilde{T}_0 & \tilde S \\
0 & \tilde{T}_1 \\
\end{pmatrix}$ be two operators in $\mathcal FB_2(\Omega)$. Also let $t_1$ and
$\tilde{t}_1$ be non-zero sections of the holomorphic Hermitian vector
bundles $E_{T_1}$ and $E_{\tilde{T}_1}$ respectively. The operators $T$ and
$\tilde{T}$ are equivalent if and only if
$\mathcal{K}_{T_0}=\mathcal{K}_{\tilde{T}_0}$ (or
$\mathcal{K}_{T_1}=\mathcal{K}_{\tilde{T}_1}$) and
$\frac{\|S(t_1)\|^2}{\|t_1\|^2}= \frac{\|\tilde
S(\tilde{t}_1)\|^2}{\|\tilde{t}_1\|^2}$.
\end{theorem}  

In any decomposition $\Big (\begin{smallmatrix}T_0 & S \\
0 & T_1 \\
\end{smallmatrix} \Big ),$
of an operator $T\in \mathcal FB_2(\Omega),$ let $t_1$ be a non zero section of holomorphic Hermitian vector bundle $E_{T_1}$. 
We  assume, without loss of generality, that $S(t_1)$ is a non zero
section of $E_{T_0}$ on some open subset of $\Omega$. Following the methods of \cite[pp. 2244]{dm}, the second fundamental form  of $E_{T_0}$ in $E_T$ is easy to compute. It is the $(1,0)$-form
$\tfrac{\mathcal{K}_{T_0}(z)}{\,\,\Big(-\mathcal{K}_{T_0}(z)+\tfrac{\|t_1(z)\|^2}{\big \|S(t_1(z) )\big \|^2}\Big)^{\!\!\stackrel{}{1}/2}}d\bar{z}.$
Thus the second fundamental form of $E_{T_0}$ in $E_T$ together with the curvature of
$E_{T_0}$ is a complete set of invariants for the operator $T$. The inclusion of the line bundle $E_{T_0}$ in the vector bundle 
$E_{T}$ of rank $2$ is the flag structure of $E_T$.


It is not easy to determine which operators in $B_n(\Omega)$ are irreducible. We show that the operators in the new class $\mathcal FB_2(\Omega)$ are always irreducible.  Indeed, if we assume $S$ is invertible, then $T$ is strongly irreducible, that is, there is no non-trivial idempotent commuting with $T.$ 
  
Recall that an operator $T$ in the Cowen-Douglas class $B_n(\Omega)$,  up to unitary equivalence, is the  adjoint of the multiplication operator {$M$} on a Hilbert space {$\mathcal H$} consisting of holomorphic functions on $\Omega^*:=\{\bar{w}: w \in \Omega\}$ possessing a reproducing kernel {$K$.}  What about operators in $\mathcal FB_2(\Omega)$? A model for these operators is described below.   

Let $\gamma=(\gamma_0,\gamma_1)$ be a holomorphic frame for the vector bundle $E_T$, $T\in \mathcal FB_2(\Omega)$.  Then the operator $T$ is unitarily equivalent to the adjoint of the multiplication operator $M$ on a reproducing kernel Hilbert space
$\mathcal{H}_{\Gamma} \subseteq {\rm Hol}(\Omega^*, \mathbb C^2)$ possessing a reproducing kernel $K_{\Gamma}:\Omega^* \times \Omega^* \to \mathbb C^{2\times 2}$ of the special form that we describe explicitly now. For $z,w\in\Omega^*$, 
\begin{eqnarray*}
K_{\Gamma}(z,w)&=&
\begin{pmatrix}
  \langle \gamma_0(\bar w),\gamma_0(\bar z)\rangle &
  \langle \gamma_1(\bar w),\gamma_0(\bar z)\rangle \\
   \langle \gamma_0(\bar w),\gamma_1(\bar z)\rangle &
   \langle \gamma_1(\bar w),\gamma_1(\bar z)\rangle \\
\end{pmatrix}\\
&=& \begin{pmatrix}
  \langle \gamma_0(\bar w),\gamma_0(\bar z)\rangle &
  \frac{\partial}{\partial \bar w}\langle \gamma_0(\bar w),\gamma_0(\bar z)\rangle \\
  \frac{\partial}{\partial z} \langle \gamma_0(\bar w),\gamma_0(\bar z)\rangle
  &
 \frac{\partial^2}{\partial z\partial \bar w} \langle \gamma_0(\bar w),
 \gamma_0(\bar z)\rangle+\langle t_1(\bar w),t_1(\bar z)\rangle \\
\end{pmatrix}, 
\end{eqnarray*}
where  $t_1$ and $\gamma_0:=S(t_1)$ are
frames of the line bundles $E_{T_1}$ and $E_{T_0}$ respectively. 
It follows that  $\gamma_1(w):=\tfrac{\partial}{\partial w}\gamma_0(w)-t_1(w)$ and that 
$t_1(w)$ is orthogonal to $\gamma_0(w)$, $w\in\Omega$.

Setting $K_0(z,w)=\langle \gamma_0(\bar w),\gamma_0(\bar z)\rangle$
and $K_1(z,w)= \langle t_1(\bar w),t_1(\bar z)\rangle$, we see that the
reproducing kernel $K_{\Gamma}$ has the form:
\begin{equation} \label{main}
K_{\Gamma}(z,w)= 
\begin{pmatrix}
  {K_0}(z,w) & \frac{\partial}{\partial \bar w}{K_0}(z,w) \\
  \frac{\partial}{\partial z}{K_0}(z,w) & \frac{\partial^2}{\partial z\partial \bar w}{K_0}
  (z,w)+{K_1}(z,w) \\
\end{pmatrix}.
\end{equation}

This special form of the kernel $K_\Gamma$ for an operator in the class $\mathcal FB_2(\Omega)$ entails that a change of frame between any two frames
$\{\gamma_0,\gamma_1\}$ and $\{\sigma_0,\sigma_1\}$ of vector
bundle $E_{T}$, which has property $\gamma_0\perp \partial
\gamma_0-\gamma_1$ and
 $\sigma_0\perp \partial \sigma_0-\sigma_1$, must be induced by a holomorphic
 $\Phi:\Omega\to \mathbb{C}^{2\times 2}$ of the form
  $\Phi=\Big ( \begin{smallmatrix}\phi & \phi{\prime} \\
0 & \phi \\
\end{smallmatrix}\Big )$ for some holomorphic function $\phi :\Omega \to \mathbb{C}.$
As an immediate corollary, we see that an unitary operator intertwining two of these operators,  represented in the form  $T:=\Big ( \begin{smallmatrix} T_0 & S \\ 0 & T_1  \end{smallmatrix}\Big )$ and $\tilde{T}:=\Big ( \begin{smallmatrix} \tilde{T}_0 & \tilde{S} \\ 0 & \tilde{T}_1  \end{smallmatrix}\Big )$, must be diagonal with respect to the implicit decomposition of the two Hilbert spaces $\mathcal H$ and $\tilde{\mathcal H}.$  As a second corollary, we see that if $T_0=\tilde{T}_0$ and $T_1=\tilde{T}_1,$ then the operators $T$ and $\tilde{T}$ are unitarily equivalent if and only if $\tilde{S} = e^{i\theta} S$ for some real $\theta.$
  
We now give examples of natural classes of operators that belong to $\mathcal FB_2(\Omega)$.  Indeed, we were led to the definition of this new class  $\mathcal F B_2(\Omega)$ of operators by trying to understand these examples better.  

An operator $T$ is called \emph{homogeneous} if $\phi(T)$ is unitarily equivalent to
$T$ for all $\phi$ in M\"{o}b which are analytic on the spectrum of
$T$.
  

If an operator $T$ is in ${B}_1(\mathbb D)$,  then $T$ is  homogeneous if and only if $\mathcal{K}_T(w)=-\lambda
(1-|w|^2)^{-2},$ for some $\lambda>0$. A model for all homogeneous operators in $B_n(\mathbb D)$ is given in \cite{AK}. We describe them for $n=2$. For $\lambda > 1$ and $\mu > 0$, set $K_0(z,w) = (1-z\bar{w})^{-\lambda}$ and $K_1(z,w) = \mu (1-z\bar{w})^{-\lambda-2}$.  An irreducible operator $T$ in $B_2(\mathbb D)$ is homogeneous if and only if it is unitarily equivalent to the adjoint of the multiplication operator on the Hilbert space $\mathcal H\subseteq {\rm Hol}(\mathbb D, \mathbb C^2)$ determined by the positive definite kernel given in equation  \eqref{main}.  The similarity as well as a unitary classification of homogeneous operators in $B_n(\mathbb D)$ were obtained in \cite{AK}  using non-trivial results from representation theory of semi-simple Lie  group.  For $n=2$, this classification is a consequence of Theorem \ref{maint}.  

An operator $T$ in $B_1(\Omega)$ acting on a  Hilbert space 
$\mathcal H$ makes it a  module over the polynomial ring via the usual point-wise multiplication.  An important tool in the study 
of these modules is the module tensor product (or, localization). This is the Hilbert module 
$J \mathcal H^{(k)}_{\rm loc}$ corresponding to the spectral sheaf 
$J \mathcal H \otimes_\mathcal P \mathbb C^k_w$, where 
$\mathcal P$ is the polynomial ring and  
\begin{enumerate}
\item $J:\mathcal H \to {\rm Hol}(\Omega, \mathbb C^k)$ is the jet map, namely, $Jf = \sum_{\ell=0}^{k-1} \partial^\ell f \otimes \varepsilon_{\ell+1},$ 
$\varepsilon_1, \ldots ,\varepsilon_{k}$ are the standard unit vectors in $\mathbb C^k$.
\item $\mathbb C^k_w$ is a $k$ - dimensional module over the polynomial ring, 
\item the module action on $\mathbb C^k_w$ is via the map $\mathcal J(w),$ see \cite[(2.8) pp. 376]{dmv}:
\renewcommand\arraystretch{.675}
$$(\mathcal Jf)(w)=\left(
                         \begin{array}{cccc}
                           f(w) & 0 & \cdots & 0 \\
                           \binom{2}{1}\partial f(w) & f(w) & \cdots & 0 \\
                           \vdots & \vdots & \ddots & \vdots \\
                           \binom{k}{1}\partial^{k-1}f(w) & \binom{k-1}{1}\partial^{k-2}f(w) & \cdots & f(w) \\
                         \end{array}
                       \right),$$
\renewcommand\arraystretch{1}                       
\hspace{-6ex} that is,  $(f,v) \mapsto (\mathcal Jf)(w) v,$ $f\in \mathcal P,\, v\in \mathbb C^k.$
\end{enumerate}
  
We now consider the localization with $k=2$. If we assume that the operator $T$ has been realized as the adjoint of the multiplication operator on a Hilbert space of holomorphc function possessing a kernel function, say $K$, then the kernel $JK^{(2)}_{\rm loc}$ for the localization (of rank $2$) given in \cite[(4.2) pp.  393]{dmv} coincides with $K_\Gamma$ of equation \eqref{main}. In this case,  we have $K_1=K=K_0$. 
The operator {$T,$} in this case, has the form {$\Big (\begin{smallmatrix} T & \binom{2}{1}\,I\\ 0 & T\\
  \end{smallmatrix}\Big ).$}
As is to be expected, using the complete set of unitary invariants given in Theorem \ref{maint}, we see that the unitary equivalence class of the Hilbert module $\mathcal H$ is in one to one correspondence with that of $J \mathcal H^{(2)}_{\rm loc}$.

Thus the class $\mathcal FB_2(\Omega)$ contains two very interesting classes of operators. For $n > 2$, we find that there are competing definitions. One of these contains the homogeneous operators and the other contains the Hilbert modules obtained from the localization. At this point, we note that most of what is said for the class $\mathcal F B_2(\Omega)$ remains valid if we assume $T_0$ is in $B_{n}(\Omega)$, $n>1$, instead of $B_1(\Omega)$. Although, now we must assume that the operator $S$ has dense range, merely assuming that it is non-zero is not enough. Also, it is no longer possible to describe a complete set of invariants for such an operator as in Theorem \ref{maint}. We proceed slightly differently to ensure a better understanding.

\begin{e-definition}Let $\mathcal FB_n(\Omega)$ be the set of all operators $T$ in the Cowen-Douglas class $B_n(\Omega)$ for which we can find operators $T_0,T_1, \ldots , T_{n-1}$ in $B_1(\Omega)$ and a decomposition of the form  
\renewcommand\arraystretch{.625}
$$T=\begin{pmatrix}
T_{0} & S_{0\,1} &S_{0\,2}&\ldots&S_{0\,n-1}\\
0&T_{1}&S_{1\,2}&\ldots&S_{1\,n-1} \\
\vdots&\ddots&\ddots&\ddots&\vdots\\
0&\ldots&0&T_{n-2}&S_{n-2\,n-1}\\
0&\ldots&\ldots&0&T_{n-1}\\
\end{pmatrix}$$
\renewcommand\arraystretch{1}
such that none of the operators $S_{i\,i+1}$ are zero and $T_i S_{i\,i+1} = S_{i\,i+1}T_{i+1}.$
\end{e-definition}

If there exists a invertible bounded linear operator $X$  intertwining any two operators, say  $T,\, \tilde{T}$ in $\mathcal FB_n(\Omega)$ ($XT = \tilde{T}X$), then we prove that $X$ must be upper triangular with respect to the decomposition mandated in the definition of the class $\mathcal FB_n(\Omega).$ 
It then follows that any unitary operator intertwining these two operators 
must be diagonal.  Thus we see that they are unitarily equivalent if and only there exists unitary operators $U_i: \mathcal H_i \to \tilde{\mathcal H}_i$, $i=0,1,\cdots n-1,$ such that $U_i^*\widetilde{T}_iU_i=T_i$ and $U_iS_{i,j}=\widetilde{S}_{i,j}U_j,$ $i<j$. The first of these conditions immediately translates into a condition on the curvature of the line bundles $E_{T_i}.$ The second condition is somewhat more mysterious and is related to a finite number of second fundamental forms inherent in our description of the operator $T.$  In what follows, we make this a little more explicit after making some additional assumptions. 

Let $T$ be an operator acting on a Hilbert space $\mathcal H.$   
Assume that there exists a representation of the form  
{
\renewcommand\arraystretch{.625}
\begin{equation}\label{decomp}
T= 
\begin{pmatrix}
T_{0}&S_{0\,1}&0&\ldots&0\\
0&T_{1}&S_{1\,2}&\ldots&0\\
\vdots&\ddots&\ddots&\ddots&\vdots \\
0&\ldots&0&T_{n-2}&S_{{n-2\,n-1}}\\
0&\ldots&0&0&T_{n-1} 
\end{pmatrix}
\end{equation}\renewcommand\arraystretch{1}}
for the operator $T$ with respect to some orthogonal decomposition $\mathcal H:=\mathcal{H}_0\oplus\mathcal{H}_1\oplus \cdots \oplus \mathcal{H}_{n-1}.$
Suppose also that the operator $T_{i}$ is in $B_1(\Omega),$ $0\leq i \leq n-1,$ the operator $S_{i-1,i}$ is non zero and $T_{k-1}S_{k-1,k}=S_{k-1,k}T_{k},$ $1\leq i\leq n-1.$ Then we show that the operator $T$ must be in the Cowen-Douglas class $B_n (\Omega).$ We can also relate the frame of the vector bundle $E_T$ to those  of the line bundles $E_{T_i},$ $i=0,1,\ldots , n-1.$ Indeed, we show that there is a frame $\{\gamma_0,\gamma_1,\cdots,\gamma_{n-1}\}$ of $E_{T}$ such that
$$t_k(w):=\gamma_k(w)+\cdots+(-1)^j\binom{k}
{i}\gamma_{k-j}^{(j)}(w)+\cdots+(-1)^{k}\gamma_0^{(k)}(w)$$ is
a non-zero section of the line bundle $E_{T_k}$ and it is orthogonal to $\gamma_i(w)$, $i=0,1,2,\ldots, k-1.$  We also have $t_{i-1}:=S_{i-1\,i}(t_i),$ $1 \leq i \leq n-1$.
In this special case, we can extract a complete set of  invariants explicitly. 
\setcounter{theorem}{1}
\begin{theorem} Pick two operators $T$ and $\tilde{T}$ which admit a decomposition of the form given in \eqref{decomp}. Find  an orthogonal frame $\{\gamma_0,t_1,\cdots,t_{n-1}\}$ (resp. $\{\tilde{\gamma}_0,\tilde{t}_{1},\cdots,\tilde{t}_{n-1}\}$) 
for the vector bundle ${\bigoplus\limits^{n}_{i=0}E_{T_i}}$ (resp. 
${\bigoplus\limits^{n}_{i=0}E_{\tilde{T}_i}}$) as above. Then the operators $T$ and $\tilde{T}$ are unitarily
equivalent if and only if
$$\mathcal K_{T_0} =\mathcal K_{\tilde{T}_0}\mbox{ 
and }
\frac{\|S_{i-1\,i}(t_i)\|^2}{\|t_{i}\|^2}=\frac{\|\tilde{S}_{i-1\,i}(\tilde{t}_i)\|^2}{\|\tilde{t}_{i}\|^2},\,\, 1\leq i\leq n-1.$$
\end{theorem}
The operator corresponding to the module action $J\mathcal H^{(k)}_{\rm loc},$ or for that matter $\mathcal H^{(k)}_{\rm loc},$ has exactly the form described above while the homogeneous operators (of rank $>2$) are not of this form, although  they meet the requirements set forth in Definition 2.   

}

\end{document}